\documentclass[11pt]{article}
\usepackage{amsfonts,amsmath,amsthm,amssymb}
\theoremstyle{definition}
\newtheorem{theorem}{Theorem}
\newtheorem{lemma}[theorem]{Lemma}
\newtheorem{cor}[theorem]{Corollary}

\newtheorem{prop}[theorem]{Proposition}

\newcommand\shc[3]{\ensuremath{{#1}_{#2}^{(#3)}}}
\newcommand\E{\mathbb{E}}

\newcommand\comment[1]{}

\title{Density of normal binary covering codes}
\author{Robert B.\ Ellis\thanks{\small
Research supported in part by NSF grant DMS-9977354.}\\
Texas A\&M University
\date{December 31, 2003} \\}

\allowdisplaybreaks

\begin{document}
\maketitle

\begin{abstract}
A binary code with covering radius $R$ is a subset $\mathcal{C}$
of the hypercube $Q_n=\{0,1\}^n$ such that every $x\in Q_n$ is
within Hamming distance $R$ of some codeword $c\in \mathcal{C}$,
where $R$ is as small as possible.  For a fixed coordinate
$i\in[n]$, define $\mathcal{C}^{(i)}_b$, for $b\in\{0,1\}$, to be
the set of codewords with a $b$ in the $i$th position.  Then
$\mathcal{C}$ is {\em normal} if there exists an $i\in[n]$ such
that for any $v\in Q_n$, the sum of the Hamming distances from $v$
to $\mathcal{C}^{(i)}_0$ and $\mathcal{C}^{(i)}_1$ is at most
$2R+1$. We newly define what it means for an asymmetric covering
code to be normal, and consider the worst case asymptotic
densities $\nu^*(R)$ and $\nu^*_+(R)$ of constant radius $R$
symmetric and asymmetric normal covering codes, respectively.
Using a probabilistic deletion method, and analysis adapted from
previous work by Krivelevich, Sudakov, and Vu, we show that both
are bounded above by $e(R\log R + \log R + \log\log R+4)$, giving
evidence that minimum size constant radius covering codes could
still be normal.
\end{abstract}

\section{Introduction}
The problem of finding a small set of $n$-bit binary string {\em
codewords} such that every $n$-bit binary string is within $R$
bit-flips of a codeword is the classical coding theory question of
finding binary {\em covering codes} of length $n$ and radius $R$.
Much effort has been made to determine the minimum or {\em
optimal} size of the smallest binary covering codes for various
values of $n$ and $R$, as well as for constant $R$ as $n$ tends to
infinity (cf.\ Chapter 12 of \cite{CHLL97}), with asymptotically
tight bounds having been achieved only in the case of $R=1$. One
method by Graham and Sloane \cite{GS85}, which has produced
best-known upper bounds on the optimal size of covering codes for
many values of $n$ and $R$, involves considering a special class
of so-called {\em normal} codes (cf.\ entries marked with ``Q'' in
Table 6.1 of \cite{CHLL97}). These codes admit to an efficient
concatenation operation, called {\em amalgamated direct sum}
(ADS), by which good longer codes are constructed from shorter
codes. In this paper, we extend this concatenation operation to
give an asymptotic upper bound on the optimal size of constant
radius normal covering codes which nearly approaches the
corresponding best-known bound for unrestricted codes. Our
extension employs a probabilistic deletion method, and a recursive
construction motivated by \cite{CEK02} and \cite{KSV03}, from
which several analytical techniques are also borrowed.  This
result provides positive evidence for an unsolved conjecture: for
general $n$ and $R$, does there exist an optimal code which is
also normal?  We also newly define normality for {\em asymmetric
codes}, in which every $n$-bit string must be obtainable from a
codeword by flipping at most $R$ 1's to 0's, and we adapt the
above-mentioned extended concatenation operation to give an
asymptotic bound on the optimal size of normal asymmetric codes
for constant $R$.

\section{Definitions and the ASDS construction}

Let $Q_n:=\{x=(x_1,x_2,\ldots,x_n):x_i\in\{0,1\}\}$ be the set of
$n$-bit strings, or binary $n$-vectors, with algebraic structure
inherited from the vector space $\mathbb{F}_2^n$ and partial
ordering inherited from the boolean lattice (i.e., $x\preceq y$
provided $x_i\leq y_i$ for all $1\leq i\leq n$).  Define the {\em
weight}, or {\em level}, of $x\in Q_n$ to be $w(x):=\sum_{i=1}^n
x_i$, that is, the number of 1's in $x$.  Define the Hamming
distance between $x$ and $y$ to be $d(x,y):=w(x-y)$; for a set
$Y\subseteq Q_n$, $d(x,Y):=\min\{d(x,y):y\in Y\}$, with
$d(x,Y)=\infty$ when $Y=\emptyset$. The {\em undirected ball} in
$Q_n$ with center $x$ and radius $R$, denoted by $B_n(x,R)$, is
the set $\{y\in Q_n:d(x,y)\leq R\}$.  We sometimes refer to such a
ball as an {\em $R$-ball}. The size of $B_n(x,R)$ is independent
of $x$ and is denoted by $b_n(R)$. The {\em covering radius} of a
set $\mathcal{C}\in Q_n$ is the smallest integer $R\geq 0$ such
that $Q_n = \cup_{c\in\mathcal{C}}B_n(c,R)$.  The usual definition
of a binary covering code, which for our purposes we refer to as a
{\em symmetric binary covering code} of length $n$ and radius $R$,
or more simply an $(n,R)$-code, is a set of {\em codewords}
$\mathcal{C}\subseteq Q_n$ with covering radius $R$. We use
$K(n,R)$ to denote the minimum size of any $(n,R)$-code. A lower
bound for $K(n,R)$ is obtained by considering that the minimum
conceivable number of $R$-balls needed to cover $Q_n$ is
$2^n/b_n(R)$, which gives the (folkloric) sphere bound
$$ K(n,R) \ \geq \ \frac{2^n}{b_n(R)} \ = \frac{2^n}{\binom{n}{\leq R}},  $$
where we define $\binom{n}{\leq R}:=\sum_{i=0}^R \binom{n}{i}$.
The sphere bound motivates the definition of the {\em density} of
an $(n,R)$-code $\mathcal{C}$, which is
$\frac{|\mathcal{C}|}{2^n/\binom{n}{\leq R}}$.  The {\em optimal
density} of an $(n,R)$-code is
$\mu(n,R):=\frac{K(n,R)}{2^n/\binom{n}{\leq R}}$, and the {\em
asymptotic worst-case density} of an $(n,R)$-code is
$$ \mu^*(R) \ := \limsup_{n\rightarrow\infty} \mu(n,R).  $$
It is known that $\mu^*(1)=1$ by Theorem 12.4.11 of \cite{CHLL97}
due to Kabatyanskii and Panchenko; whether $\mu^*(R)=1$ for
constant $R\neq 1$ is a central conjecture in coding theory.

In order to define asymmetric covering codes, we first define
upward and downward directed $R$-balls.  An {\em upward directed
ball} in $Q_n$ with center $x$ and radius $R$ is defined as the
set $B_n^+(x,R):=B_n(x,R)\cap \{y\in Q_n:x\preceq y\}$, and the
corresponding {\em downward directed ball} is
$B_n^-(x,R):=B_n(x,R)\cap\{y\in Q_n:y\preceq x\}$. We write
$b_n^+(x,R)$ or $b_n^-(x,R)$ for the sizes of the upward or
downward directed $R$-balls centered at $x\in Q_n$, respectively,
and sometimes instead write $b_n^+(l,R)$ or $b_n^-(l,R)$, where
$l$ is the weight $w(x)$ of $x$, since directed ball size depends
only on $n$, $R$, and the weight of the center $x$.  In
particular,
$$ b_n^+(l,R) \ = \ b_n^-(n-l,R) \ = \
\binom{n-l}{\leq R}.  $$ The {\em asymmetric distance} $d^+(x,Y)$
between a vector $x\in Q_n$ and a set $Y\subseteq Q_n$ is defined
by $d^+(x,Y):=\min\{d(x,y):y\in Y\mbox{ and }x\preceq y\}$, to
reflect the fact that $x$ can be covered by $B_n^-(y,R)$ provided
that $d^+(x,y)\leq R$. A set $\mathcal{C}\subseteq Q_n$ {\em
downward $R$-covers} $Q_n$ provided that $Q_n =
\cup_{c\in\mathcal{C}} B_n^-(c,R)$, and the asymmetric covering
radius of $\mathcal{C}$ is the smallest $R$ for which
$\mathcal{C}$ downward $R$-covers $Q_n$.  We say that such a set
$\mathcal{C}$ with asymmetric covering radius $R$ is an {\em
asymmetric binary covering code} of length $n$ and radius $R$, or
more simply, an $(n,R)^+$-code. Analogous to the notation for
symmetric codes, we define $K^+(n,R)$ to be the minimum size of an
$(n,R)^+$-code. Since the typical downward directed $R$-ball size
in $Q_n$ is $\binom{\lfloor n/2\rfloor}{\leq R}$, following
\cite{KSV03} we define the density of an $(n,R)^+$-code
$\mathcal{C}$ to be $\frac{|\mathcal{C}|}{2^n/\binom{\lfloor
n/2\rfloor}{\leq R}}$; an alternate definition for small values of
$n$ and $R$ is given in Theorem 2 of \cite{CEK02}. The optimal
density of an $(n,R)^+$-code is
$\mu_+(n,R):=\frac{K^+(n,R)}{2^n/\binom{\lfloor n/2\rfloor}{\leq
R}}$, and the asymptotic worst-case density of an $(n,R)^+$-code
is
$$ \mu_+^*(R) \ := \ \limsup_{n\rightarrow\infty} \mu_+(n,R). $$
For properties of $(n,R)^+$-codes, especially for constant $R$ or
constant $n-R$, see \cite{CEK02}.

The {\em concatenation} of two vectors $x\in Q_n$ and $y\in
Q_{n'}$ is the vector $(x,y)\in Q_{n+n'}$ determined by
$(x,y):=(x_1,\ldots,x_n,y_1,\ldots,y_{n'})$.  The {\em direct sum}
of two sets $X\subseteq Q_n$ and $Y\subseteq Q_{n'}$ is $X\oplus
Y:=\{(x,y):x\in X,y\in Y\}\subseteq Q_{n+n'}$. The following
proposition is straightforward and presented without proof, as it
is well-known in the symmetric case.
\begin{prop}[Direct sum of codes]
Let $\mathcal{C}$ be an $(n,R)$-code ($(n,R)^+$-code), and let
$\mathcal{C}'$ be an $(n',R')$-code ($(n',R')^+$-code). Then
$\mathcal{C}\oplus \mathcal{C'}$ is an $(n+n',R+R')$-code
($(n+n',R+R')^+$-code).
\end{prop}
We have reminded the reader of the direct sum construction because
it is the basis of the amalgamated direct sum and amalgamated
semi-direct sum constructions to be defined.

\subsection{Normal codes}

We now present normal symmetric covering codes, introduced in
\cite{GS85}; our notation follows that of Chapter 4 in
\cite{CHLL97}. Let $[n]:=\{1,\ldots,n\}$.  For a fixed coordinate
$i\in[n]$ and a set $X\subseteq Q_n$, define
$\shc{X}{0}{i}:=\{x\in X:x_i=0\}$, and $\shc{X}{1}{i}:=\{x\in
X:x_i=1\}$; thus $\shc{\mathcal C}{0}{i}$ and $\shc{\mathcal
C}{1}{i}$ partition a code $\mathcal{C}\subseteq Q_n$ based on the
$i$th codeword coordinate. The norm of $\mathcal{C}$ with respect
to the $i$th coordinate is
$$ N^{(i)}:= \max_{x\in Q_n}
    \left\{d(x,\shc{\mathcal C}{0}{i})
    +d(x,\shc{\mathcal C}{1}{i})\right\}.  $$
The {\em minimum norm} of a code $\mathcal{C}$ with length $n$ is
defined to be
$$ N_{\min}(\mathcal{C}) := \min_{i\in[n]} N^{(i)}.$$
A code $\mathcal{C}$ has {\em norm} $N$ provided
$N_{\min}(\mathcal{C})\leq N$. In other words, $\mathcal{C}$ has
norm $N$ provided there is a coordinate $i$ such that
$d(x,\shc{\mathcal C}{0}{i}) +d(x,\shc{\mathcal C}{1}{i})\leq N$
for all $x\in Q_n$.  A code with covering radius $R$ is {\em
normal} provided it has norm $N=2R+1$ and its minimum norm
$N_{\min}$ is $2R+1$ or $2R$, since if a code has norm $N$, its
covering radius is $R\leq N/2$. If $N^{(i)}\leq 2R+1$, then
coordinate $i$ is {\em acceptable} with respect to $2R+1$. We
shall refer to such a code as a {\em symmetric normal}
$(n,R)$-code, or equivalently a normal $(n,R)$-code. Define
$K_{\nu}(n,R)$ to be the size of the smallest normal $(n,R)$-code,
$\nu(n,R):=\frac{K_{\nu}(n,R)}{2^n/\binom{n}{\leq R}}$ to be the
optimal density of a normal $(n,R)$-code, and
$\nu^*(R):=\limsup_{n\rightarrow\infty}\nu(n,R)$ to be the
asymptotic worst-case density of a normal $(n,R)$-code. By Theorem
4.4.2 of \cite{CHLL97} due to Honkala and H{\"a}m{\"a}l{\"a}inen,
and independently van Wee \cite{vW90}, all optimal $(n,1)$-codes
with length $n\geq 3$ are normal. Therefore $\nu^*(1)=\mu^*(1)=1$,
but it is unknown whether equality holds for $R>1$.

The {\em asymmetric norm} of a code is newly defined here and is
similar to the (symmetric) norm above. Notation which coincides
with that of the symmetric norm will be made clear from context.
 The asymmetric norm of a code $\mathcal{C}$ of length $n$ with
respect to coordinate $i$ is
$$ N^{(i)}(\mathcal{C}) :=
\max\left\{\max_{x\in \shc{(Q_{n})}{0}{i}}
\left\{d^+(x,\shc{\mathcal C}{0}{i})
    + d^+(x,\shc{\mathcal C}{1}{i})\right\}, \max_{x\in \shc{(Q_n)}{1}{i}}
    \left\{2\cdot d^+(x,\shc{\mathcal C}{1}{i})+1\right\} \right\}; $$
The departure from the definition of the (symmetric) norm with
respect to coordinate $i$ is due to the fact that a vector
$x\in\shc{(Q_n)}{1}{i}$ cannot be covered by any downward directed
ball centered in $\shc{(Q_n)}{0}{i}$. The {\em minimum asymmetric
norm} $N_{\min}$ of $\mathcal{C}$ is
$$ N_{\min}(\mathcal{C}) := \min_{i\in[n]} N^{(i)}(\mathcal{C}).  $$
Therefore if a code $\mathcal{C}$ has asymmetric norm $N$, there
is a coordinate $i$ such that all words $x$ with $x_i=0$ satisfy
$d^+(x,\shc{\mathcal C}{0}{i})+d^+(x,\shc{\mathcal C}{1}{i})\leq
N$, and all words $x$ with $x_i=1$, for which $d^+(x,\shc{\mathcal
C}{0}{i})=\infty$, satisfy $d^+(x,\mathcal{C}_1^{(i)})\leq
(N-1)/2$. An $(n,R)^+$-code is {\em asymmetric normal}, or simply
{\em normal} if the context is clear, provided it has asymmetric
norm  $N=2R+1$ and its minimum asymmetric norm $N_{\min}$ is
$2R+1$ or $2R$. If $N^{(i)}\leq 2R+1$, then coordinate $i$ is {\em
acceptable} with respect to $2R+1$. Define $K^+_{\nu}(n,R)$ to be
the size of the smallest normal $(n,R)^+$-code,
$\nu_+(n,R):=\frac{K^+_{\nu}(n,R)}{2^n/\binom{\lfloor
n/2\rfloor}{\leq R}}$ to be the optimal density of a normal
$(n,R)^+$-code, and
$\nu_+^*(R):=\limsup_{n\rightarrow\infty}\nu_+(n,R)$ to be the
asymptotic worst-case density of a normal $(n,R)^+$-code.

\subsection{Amalgamated direct sum (ADS) of normal codes}

Two normal codes can be concatenated in a more efficient
construction than the basic direct sum.  The construction is the
same regardless of whether considering symmetric or asymmetric
codes, and so we present the two cases simultaneously in the
following theorem, the symmetric case of which is due to Graham
and Sloane \cite{GS85}. The theorem in the symmetric case is often
stated in terms of the covering radius, but in this paper the norm
is more central to our purpose.

\begin{theorem}[ADS of normal codes]\label{thm:ADS}
Let $A$ be a normal symmetric (asymmetric) code of length $n_A$
and norm $N_A$ with the last coordinate acceptable, and let $B$ be
a normal symmetric (asymmetric) code of length $n_B$ and norm
$N_B$ with the first coordinate acceptable.  Then their {\em
amalgamated direct sum} (ADS)
$$ A \dot\oplus B \ := \ \{(a,0,b):(a,0)\in A,(0,b)\in B\}
     \cup \{(a,1,b):(a,1)\in A, (1,b)\in B\}
     $$
is a normal symmetric (asymmetric) code of length $n_A+n_B-1$ and
norm $N_A+N_B-1$ with respect to coordinate $n_A$.
\end{theorem}
\begin{proof}
The proof of the symmetric case essentially appears in the proof
of Theorem 4.1.8 and the remarks following Theorem 4.1.14, both of
\cite{CHLL97}. We now adapt the same proof for the asymmetric
case, from which the reader may easily reconstruct the symmetric
case.

Let $\mathcal{C}=A\dot\oplus B$.  Then $\mathcal{C}$ clearly has
length $n_A+n_B-1$, as it is constructed by overlapping a single
coordinate of $A$ and $B$. Let $z\in \mathcal{C}$.  First suppose
$z=(x,0,y)$, where $(x,0)\in Q_{n_A}$.  Computing, we have
\begin{eqnarray}
d^+(z,\shc{\mathcal C}{0}{n_A})+d^+(z,\shc{\mathcal C}{1}{n_A})
    & \leq & d^+((x,0),\shc{A}{0}{n_A}) + d^+((0,y),\shc{B}{0}{1})
    \nonumber \\
& & + \ d^+((x,0),\shc{A}{1}{n_A}) +
    d^+((0,y),\shc{B}{1}{1})-1 \nonumber \\
& \leq & N_A + N_B - 1. \nonumber
\end{eqnarray}
Now suppose $z=(x,1,y)$, where $(x,1)\in Q_{n_A}$. Then we have
\begin{eqnarray}
2\cdot d^+(z,\shc{\mathcal C}{1}{n_A})+1 & \leq &
    \big(2\cdot d^+((x,1),\shc{A}{1}{n_A})+1\big) +
    \big(2\cdot d^+((1,y),\shc{B}{1}{1})+1\big)-1
    \nonumber\\
& \leq & N_A+N_B-1.\nonumber
\end{eqnarray}
Therefore $\mathcal{C}$ has asymmetric norm $N_A+N_B-1$ with
respect to coordinate $n_A$.
\end{proof}
The size of $A\dot\oplus B$ depends on the relative sizes of
$\shc{A}{0}{n_A}$ versus $\shc{A}{1}{n_A}$ and of $\shc{B}{0}{1}$
versus $\shc{B}{1}{1}$. We define a code $\mathcal{C}$ to be {\em
balanced} if
$|\shc{\mathcal{C}}{0}{i}|=|\shc{\mathcal{C}}{1}{i}|$, where $i$
is the coordinate with respect to which the ADS is taken. The
major consequence of Theorem \ref{thm:ADS} for code density is as
follows. Two codes $A$ and $B$ of lengths $n_A$ and $n_B$ and
covering radii $R_A$ and $R_B$, respectively, form a direct sum of
size $|A|\cdot |B|$, length $n_A+n_B$, and covering radius
$R_A+R_B$. If in addition both codes are normal and at least one
is balanced, their amalgamated direct sum is of size $|A|\cdot
|B|/2$, length $n_A+n_B-1$, and covering radius at most $R_A+R_B$.
Since
$$
|A|\cdot |B|\frac{\binom{n}{\leq R}}{2^n} \ > \
    \frac{|A|\cdot |B|}{2}\frac{\binom{n-1}{\leq
    R}}{2^{n-1}}\qquad \mbox{and} \qquad
|A|\cdot |B|\frac{\binom{\lfloor n/2\rfloor}{\leq R}}{2^n}
    \ \geq \ \frac{|A|\cdot |B|}{2}\frac{\binom{\lfloor (n-1)/2
    \rfloor}{\leq R}}{2^{n-1}},
$$
the density of the direct sum code is at least as large as that of
the corresponding ADS code in both the symmetric and asymmetric
case.

\subsection{Amalgamated semi-direct sum (ASDS) of normal codes}

We now define the central construction of this paper, the
amalgamated semi-direct sum.  The idea behind this construction is
as follows.  With length $n$ fixed, and target norm $N$ (and
implicitly radius $R\leq N/2$), we probabilistically choose a
candidate code $S$. Any strings $x\in Q_n$ which violate the
target norm $N$ in coordinate $n$ contribute to a ``patch'' $T$.
Together, this ``patched'' code $(S,T)$ can be incorporated into a
modified amalgamated direct sum resulting in a longer code with
some desired norm, which in turn bounds the covering radius of the
resulting code.

More formally, for a fixed $N>0$, a {\em norm $N$-patched
symmetric code} of length $n$ is a 2-tuple $(S,T)$, where
$S,T\subseteq Q_n$, such that there exists a coordinate $i\in[n]$
so that for all $x\in Q_n$ either
\begin{itemize}
\item[(I)] $d(x,S_0^{(i)})+d(x,S_1^{(i)})\leq N$, or \item[(II)]
$\{x,x+e_i\}\subseteq T$, where $x+e_i$ is $x$ with the
    $i$th coordinate flipped.
\end{itemize}
When $N$ and $n$ are clear from context, the terminology {\em
norm-patched code} may also be used. Any coordinate $i$ achieving
these properties is called {\em acceptable} for $(S,T)$ with
respect to $N$. If a vector $v\in Q_n$ violates condition (I), we
say it is {\em missed} by $S$ with respect to coordinate $i$. Note
that if $(S,T)$ is a norm $N$-patched code, then $S\cup T$ is a
normal $(n,R)$-code with radius $R\leq \lfloor N/2\rfloor$.

A {\em norm $N$-patched asymmetric code} of length $n$ is defined
similarly, except that $(S,T)$ must satisfy for some coordinate
$i\in[n]$ the following altered conditions: for all $x\in
\shc{(Q_n)}{0}{i}$, either
$d^+(x,\shc{S}{0}{i})+d^+(x,\shc{S}{1}{i})\leq N$ or
$\{x,x+e_i\}\subseteq T$; and for all $x\in \shc{(Q_n)}{1}{i}$,
either $2\cdot d^+(x,\shc{S}{1}{i})+1\leq N$ or $x\in T$. A vector
$x\in \shc{(Q_n)}{0}{i}$ is {\em missed} by $S$ with respect to
coordinate $i$ provided
$d^+(x,\shc{S}{0}{i})+d^+(x,\shc{S}{1}{i})>N$, and a vector $x\in
\shc{(Q_n)}{1}{i}$ is {\em missed} by $S$ w.r.t.\ $i$ provided
$2\cdot d^+(x,\shc{S}{1}{i})+1> N$. With these definitions we have
the following new theorem.

\begin{theorem}[ASDS of norm-patched and normal codes]\label{thm:ASDS}
Suppose $(S,T)$ is a norm $N$-patched symmetric (asymmetric) code
of length $n$ with coordinate $n$ acceptable, $K_1$ is a symmetric
(asymmetric) code of length $n'$ and norm $N'$ with first
coordinate acceptable, and $K_2$ is a symmetric (asymmetric) code
of length $n'$ and norm $N+N'-1$ with first coordinate acceptable.
Then the amalgamated semi-direct sum
\begin{eqnarray}
(S,T)\dot\boxplus (K_1,K_2) \ := \ \left(S\dot\oplus
K_1\right)\cup \left(T\dot\oplus K_2\right) \nonumber
\end{eqnarray}
is a symmetric (asymmetric) code of length $n+n'-1$ and norm
$N+N'-1$ with coordinate $n$ acceptable.
\end{theorem}
\begin{proof}
First, consider the symmetric case.  Define
$\mathcal{C}=(S,T)\dot\boxplus (K_1,K_2)$ and let $z\in
Q_{n+n'-1}$. Suppose $z=(x,0,y)$ where $(x,0)\in Q_n$.  If
$d((x,0),\shc{S}{0}{n})+d((x,0),\shc{S}{1}{n})\leq N$, then we
have
\begin{eqnarray}
d(z,\shc{\mathcal{C}}{0}{n})+d(z,\shc{\mathcal{C}}{1}{n}) & \leq &
    d((x,0),\shc{S}{0}{n}) + d((0,y),\shc{(K_1)}{0}{1})
    \nonumber \\
& & + \ d((x,0),\shc{S}{1}{n}) +
    d((0,y),\shc{(K_1)}{1}{1})-1 \nonumber \\
& \leq & N + N' - 1. \nonumber
\end{eqnarray}
Otherwise we must have $\{(x,0),(x,1)\}\subseteq T$, so that
\begin{eqnarray}
d(z,\shc{\mathcal{C}}{0}{n})+d(z,\shc{\mathcal{C}}{1}{n}) & \leq &
    d((x,0),\shc{T}{0}{n}) + d((0,y),\shc{(K_2)}{0}{1})
    \nonumber \\
& & + \ d((x,0),\shc{T}{1}{n}) +
    d((0,y),\shc{(K_2)}{1}{1})-1 \nonumber \\
& \leq & 1 + (N+N'-1) -1 \ = \ N+N'-1.
    \nonumber
\end{eqnarray}
That
$d(z,\shc{\mathcal{C}}{0}{n})+d(z,\shc{\mathcal{C}}{1}{n})\leq
N+N'-1$ when $z$ is of the form $(x,1,y)$ follows by an analogous
verification, proving the theorem in the symmetric case.

For the asymmetric case, the proof that any $z$ of the form
$(x,0,y)$ for $(x,0)\in Q_n$ satisfies
$d^+(z,\shc{\mathcal{C}}{0}{n})+d^+(z,\shc{\mathcal{C}}{1}{n})\leq
N+N'-1$ is nearly identical to the symmetric case and is omitted.
Now suppose $z$ is of the form $(x,1,y)$ where $(x,1)\in Q_n$. If
$2\cdot d^+((x,1),\shc{\mathcal{C}}{1}{n})+1\leq N$, then
\begin{eqnarray}
2\cdot d^+(z,\shc{\mathcal C}{1}{n})+1 & \leq &
    \big(2\cdot d^+((x,1),\shc{S}{1}{n})+1\big) +
    \big(2\cdot d^+((1,y),\shc{(K_1)}{1}{1})+1\big)-1
    \nonumber\\
& \leq & N+N'-1.\nonumber
\end{eqnarray}
Otherwise we must have $(x,1)\in T$, so that
\begin{eqnarray}
2\cdot d^+(z,\shc{\mathcal C}{1}{n})+1 & \leq & 2\cdot
    d^+((x,1),\shc{T}{1}{n}) +
    \big(2\cdot d^+((1,y),\shc{(K_2)}{1}{1})+1\big) \nonumber\\
& \leq & N+N'-1; \nonumber
\end{eqnarray}
therefore the theorem also holds in the asymmetric case.
\end{proof}

Again, we chose to present the theorem in terms of norms of codes
rather than radii to suit our purpose in developing the main
density theorems of the next two sections.  Additionally, it will
be convenient to choose $S$ and $T$ to be balanced with respect to
the acceptable coordinate, so that the size of the resulting ASDS
can be readily determined.

\section{Asymptotic density of normal symmetric
codes\label{sec:AsymptDensity}}

We now present the main theorem on the asymptotic worst-case
density of constant radius normal symmetric codes. The framework
and analysis of the theorem borrows from that of Theorem 1.2 (and
Corollaries 1.3-1.4) of \cite{KSV03} in the following sense.  We
develop here a more careful probabilistic deletion method in Lemma
\ref{lem:delMethod} for selecting a norm-patched code $(S,T)$,
which is tailored for normal codes and our ASDS construction.  We
must also compute a preliminary asymptotic bound on the sizes of
$|S|$ and $|T|$ in Corollary \ref{cor:patchLim} before employing a
recursive ASDS construction. We then adapt Theorem 1.2 of
\cite{KSV03} and its supporting analysis from the setting of
unrestricted codes and the so-called semi-direct sum, to the case
of normal codes and our ASDS construction, in order to obtain the
main density theorem on $\nu^*(R)$. The proof of Theorem
\ref{thm:densitySym} follows these supporting results.

\begin{theorem}\label{thm:densitySym}  Let $R\geq 2$.  Then
$$ \nu^*(R) \ \leq \ e(R\log{R}+\log{R}+\log{\log{R}}+4).
$$
\end{theorem}

\begin{lemma}[Selection of a norm-patched code]
\label{lem:delMethod}
For every positive constant $x$ and positive integer $N\leq n$,
there exist (disjoint) sets
$S_0\subseteq \shc{(Q_n)}{0}{n}$ and $S_1\subseteq \shc{(Q_n)}{1}{n}$ each
of size at most
$$ \frac{x2^{n-1}}
    {b_{n-1}\left(\left\lfloor\frac{N-1}{2}\right\rfloor\right)
    +b_{n-1}\left(\left\lceil\frac{N-1}{2}\right\rceil-1\right)}
$$
and a set $T\subseteq Q_n$ of size at most $\tau(n,N,x):=$

\begin{align}
2^{n+1} &
    \sum_{i=0}^{N-1} \exp{\left( -x \frac{b_{n-1}(i-1)
    +b_{n-1}(N-i-1)}
    {b_{n-1}\left(\left\lfloor\frac{N-1}{2}\right\rfloor\right)
    +b_{n-1}\left(\left\lceil\frac{N-1}{2}\right\rceil-1\right)}
    + \frac{b_{n-1}(i-1)+b_{n-1}(N-i-1)}{2^{n-1}}
    \right)} \nonumber \\
& + \ 2^{n+1} \exp{\left(-x\frac{b_{n-1}(N-1)}
    {b_{n-1}\left(\left\lfloor\frac{N-1}{2}\right\rfloor\right)
    +b_{n-1}\left(\left\lceil\frac{N-1}{2}\right\rceil-1\right)}
    + \frac{b_{n-1}(N-1)}{2^{n-1}}    \right)},
    \label{eqn:expPatchFull}
\end{align}
such that $(S_0\cup S_1,T)$ is a balanced norm $N$-patched
symmetric code.
\end{lemma}
\begin{proof}
Let
$$ k=\left\lfloor\frac{x2^{n-1}}
    {b_{n-1}\left(\left\lfloor\frac{N-1}{2}\right\rfloor\right)
    +b_{n-1}\left(\left\lceil\frac{N-1}{2}\right\rceil-1\right)}
    \right\rfloor,
$$
and choose uniformly at random subsets $S_0\subseteq
\shc{(Q_n)}{0}{n}$ and $S_1\subseteq \shc{(Q_n)}{1}{n}$ each of
size $k$.  A vector $v\in Q_n$ is missed by $S$ if
$d(v,S_0)+d(v,S_1)>N$; otherwise, there exists an
$i\in\{0,1,\ldots,N\}$ such that $d(v,S_0)=i$ and $d(v,S_1)\leq
N-i$.
For $b\in\{0,1\}$ classify the missed vertices as follows:
\begin{eqnarray}
B_{b,1} & := & \{u\in\shc{(Q_n)}{b}{n}:0\leq d(u,S_b)< N
    \mbox{ and } d(u,S_{1-b})>N-d(u,S_b)\} \nonumber \\
B_{b,2} & := & \{u\in\shc{(Q_n)}{b}{n}:d(u,S_b)\geq N\}. \nonumber
\end{eqnarray}
Let the patch be the balanced set
$$ T \ = \ \bigcup_{b\in\{0,1\}} \left(B_{b,1}\cup B_{b,2}
    \right)+\{0,e_n\}, $$
where addition is done by taking all possible combinations of one
vector from each set and adding coordinate-wise mod 2. Thus $T$
contains all missed vertices, and $S\cup T$ is a norm $N$-patched
code. By linearity of expectation and symmetry with respect to the
$n$th coordinate,
\begin{eqnarray}
\E(|T|) & = & \sum_{b\in\{0,1\}} 2\E(|B_{b,1}|)+2\E(|B_{b,2}|)
    \ = \  4\E(|B_{0,1}|)+4\E(|B_{0,2}|) \nonumber\\
& = & 4\cdot 2^{n-1} \sum_{i=0}^{N-1}
\mbox{Pr}[d(v,S_0)=i\,|\,v\in\shc{(Q_n)}{0}{n}]
    \cdot\mbox{Pr}[d(v,S_{1})>N-i\,|\,v\in\shc{(Q_n)}{0}{n}]    \nonumber\\
& & + \ 4\cdot 2^{n-1} \mbox{Pr}[d(v,S_0)\geq
N\,|\,v\in\shc{(Q_n)}{0}{n}].\nonumber
\end{eqnarray}
For $v\in\shc{(Q_n)}{0}{n}$,
$\mbox{Pr}[d(v,S_0)=i]=\mbox{Pr}[d(v,S_0)>i-1]-\mbox{Pr}[d(v,S_0)>i]$,
and for fixed $i$, $\mbox{Pr}[d(v,S_0)>i-1]$ dominates
$\mbox{Pr}[d(v,S_0)>i]$ as $n\rightarrow \infty$; therefore we
estimate $\E(|T|)$ by
\begin{eqnarray}
\E(|T|) & \leq & 2^{n+1} \sum_{i=0}^{N-1}
\mbox{Pr}[d(v,S_0)>i-1\,|\,v\in\shc{(Q_n)}{0}{n}]
    \cdot\mbox{Pr}[d(v,S_{1})>N-i\,|\,v\in\shc{(Q_n)}{0}{n}]    \nonumber\\
& & + \ 2^{n+1} \mbox{Pr}[d(v,S_0)\geq
N\,|\,v\in\shc{(Q_n)}{0}{n}]. \label{eqn:expPatch}
\end{eqnarray}
Suppose $0\leq j< N$. For $v\in\shc{(Q_n)}{0}{n}$, if $d(v,S_{b})$
is to be more than $j$, then $S_{0}$ must not contain any of the
vertices in $B_n(v,j)\cap \shc{(Q_n)}{0}{n}$. This intersection
can be reached from $v$ by fixing the $n$th coordinate of $v$ and
changing at most $j$ of the remaining $n-1$ coordinates. In
particular, $|B_n(v,j)\cap \shc{(Q_n)}{0}{n}|=b_{n-1}(j)$. Along
with the corresponding computation for $d(v,S_{1})$, we have
\begin{eqnarray}
\mbox{Pr}[d(v,S_{0})>j\,|\,v\in\shc{(Q_n)}{0}{n}] & = &
\binom{2^{n-1}-b_{n-1}(j)}{k}\Big/
    \binom{2^{n-1}}{k}, \qquad \mbox{and}\nonumber \\
\mbox{Pr}[d(v,S_{1})>j\,|\,v\in\shc{(Q_n)}{0}{n}] & = &
\binom{2^{n-1}-b_{n-1}(j-1)}{k}
    \Big/ \binom{2^{n-1}}{k} .\nonumber
\end{eqnarray}
Now the bound on $\E(|T|)$ in (\ref{eqn:expPatch}) becomes
\begin{eqnarray}
\E(|T|) & \leq & 2^{n+1} \binom{2^{n-1}}{k}^{-2}
    \sum_{i=0}^{N-1} \binom{2^{n-1}-b_{n-1}(i-1)}{k}
    \binom{2^{n-1}-b_{n-1}(N-i-1)}{k} \nonumber \\
& & + 2^{n+1} \binom{2^{n-1}}{k}^{-1}
    \binom{2^{n-1}-b_{n-1}(N-1)}{k}.    \label{eqn:expPatchBinom}
\end{eqnarray}
Using the estimate
\begin{eqnarray}
\frac{\binom{m-d}{k}}{\binom{m}{k}} & = &
    \frac{(m-d)\cdots (m-d-k+1)}{m\cdots(m-k+1)} \nonumber \\
& \leq & \left(\frac{m-d}{m}\right)^k \ = \
    \left(1-\frac{d}{m}\right)^k
\ \leq \  e^{-k\frac{d}{m}}
        \label{eqn:binomEst}
\end{eqnarray}
borrowed from the proof of Lemma 2.2 of \cite{KSV03},
(\ref{eqn:expPatchBinom}) becomes
\begin{eqnarray}
\E(|T|) & \leq & 2^{n+1}
    \sum_{i=0}^{N-1} \exp{\left( -k \frac{b_{n-1}(i-1)
    +b_{n-1}(N-i-1)}{2^{n-1}}\right)} \nonumber \\
& & \ + \ 2^{n+1} \exp{\left(-k\frac{b_{n-1}(N-1)}{2^{n-1}}
    \right)} \nonumber \\
& \leq & \tau(n,N,x). \nonumber
\end{eqnarray}
Since there exists a $T$ of size at most $\E(|T|)$, the result
follows.
\end{proof}

In practice, what is important is the expected size of the patch
$T$ as $n\rightarrow\infty$.  We have the following asymptotic
upper bounds on $|S|$, and on $|T|$ via $\tau(n,N,x)$.

\begin{cor}\label{cor:patchLim}
Let $N\geq 2$ be fixed.  Then the asymptotic size of $S:=S_0\cup S_1$
in Lemma \ref{lem:delMethod} is
$$
|S| \ \sim \
    \left\{\begin{array}{cl}
    x \frac{2^{n}}{b_{n-1}((N-1)/2)}, & \mbox{if $N$ is odd} \\
    x\frac{2^{n-1}}{b_{n-1}(N/2-1)}, & \mbox{if $N$ is even,}
    \end{array}\right.
$$
and the size of the patch $T$ is bounded above asymptotically by
$$ \tau(n,N,x) \ \sim \
    \left\{\begin{array}{cl}
    2^{n+2}e^{-x}, & \mbox{if $N$ is odd} \\
    2^{n+1}e^{-x}, & \mbox{if $N$ is even.}
    \end{array}\right.
$$
\end{cor}
\begin{proof}
The calculation for $|S|$ is easily verified. For the size of $T$,
note that for constant $R$ the asymptotic size of an $R$-ball in
$Q_n$ is $b_n(R)\sim n^R/R!$.  The proof proceeds by identifying
which exponential terms $\exp(\cdot)$ in (\ref{eqn:expPatchFull})
are not swallowed in the limit. If $N$ is odd, then
$\left\lfloor\frac{N-1}{2}\right\rfloor= \frac{N-1}{2}$ and
$\left\lceil\frac{N-1}{2}\right\rceil-1=\frac{N-1}{2}-1$. The only
terms which survive are the
$i=\left\lfloor\frac{N-1}{2}\right\rfloor,
\left\lfloor\frac{N-1}{2}\right\rfloor+1$ terms of the summation
in (\ref{eqn:expPatchFull}), which each converge to $\exp{(-x)}$.
If $N$ is even, then $\left\lfloor\frac{N-1}{2}\right\rfloor=
\left\lceil\frac{N-1}{2}\right\rceil-1=\frac{N}{2}-1$, and the
only exponential term of the summation in (\ref{eqn:expPatchFull})
which does not vanish corresponds to $i=\frac{N}{2}$, and also
converges to $\exp{(-x)}$.  For all other exponential terms in
both cases, the numerator dominates since at least one of the two
balls has radius larger than
$\max\{\left\lfloor\frac{N-1}{2}\right\rfloor,
\left\lceil\frac{N-1}{2}\right\rceil-1\}$.
\end{proof}

The following technical lemma, due to Krivelevich, Sudakov, and Vu
\cite[Lemma 2.1]{KSV03}, allows a tight analysis of the upper
bound on $\nu^*(R)$ given by a recursive ASDS construction.  We
quote the lemma without proof and then continue to the proof of
the main theorem in the symmetric case.

\begin{lemma}[Krivelevich, Sudakov, Vu]\label{lem:KSV}
Let $(f_n)$, $(a_n)$, $(b_n)$ and $(s_n)$ be sequences of positive
numbers where
$$ \limsup_{n\rightarrow\infty} f_n\leq f, \quad
\limsup_{n\rightarrow\infty} a_n\leq a, \quad
\limsup_{n\rightarrow\infty} b_n\leq b <1  $$ and
$$ s_n\leq a_n f_{\lfloor n/y\rfloor}+b_ns_{\lfloor n/y\rfloor}$$
where $y>1$ is a constant.  Then
$$ \limsup_{n\rightarrow\infty} s_n \leq \frac{af}{1-b}.  $$
\end{lemma}

\begin{proof}[Proof of Theorem \ref{thm:densitySym}]
Let $n$ be sufficiently large ($n\geq R$ suffices), and let
$n_1=\lfloor n/R\rfloor$ and $n_1'=n-n_1+1$. The selection of
these particular parameters in the bounding of $\mu^*(R)$ is due
to \cite{KSV03}, and we find them to be suitable for the ASDS
construction as well.  We use Lemma \ref{lem:delMethod} to select
a length $n_1'$ balanced norm $(2R-1)$-patched code $(S,T)$, where
$|S|$ and $|T|$ are bounded above as given in the lemma. Let $K_1$
be an optimal normal $(n_1,1)$-code, and let $K_2$ be an optimal
normal $(n_1,R)$-code. Now perform the ASDS of $(S,T)$ with
$(K_1,K_2)$. By Theorem \ref{thm:ASDS}, the resulting code is
length $n$ and has norm $2R+1$, and so has covering radius at most
$R$. Therefore there exists a normal $(n,R)$-code with size at
most $|(S,T) \dot\boxplus(K_1,K_2)|$, and the optimal density of
such a code is
\begin{align}
\nu( & n, R) \ \leq  \
\left(\frac{|S||K_1|}{2}+\frac{|T||K_2|}{2}\right)
    \frac{\binom{n}{\leq R}}{2^n} \nonumber\\
& \leq \ \frac{1}{2}
    \frac{x2^{n_1'}}{b_{n_1'-1}(R-1)+b_{n_1'-1}(R-2)}\nu(n_1,1)
    \frac{2^{n_1}\binom{n}{\leq R}}{\binom{n_1}{\leq 1}2^n} \ + \
    \tau(n_1',2R-1,x) \frac{1}{2}\nu(n_1,R)
    \frac{2^{n_1}\binom{n}{\leq R}}{\binom{n_1}{\leq R}2^n}.
    \nonumber
\end{align}
Define $s_n:=\nu(n,R)$, $f_n:=\nu(n,1)$,
$$a_n \ := \ \frac{1}{2}\frac{x2^{n_1'}}{b_{n_1'-1}(R-1)+b_{n_1'-1}(R-2)}
    \frac{2^{n_1}\binom{n}{\leq R}}{\binom{n_1}{\leq R_1}2^n},  \quad
    \mbox{and} \quad
    b_n \ := \ \frac{1}{2}\frac{2^{n_1}\binom{n}{\leq R}}
    {\binom{n_1}{\leq R}2^n}
        \tau(n_1',2R-1,x);$$
note that
\begin{equation}
\limsup_{n\rightarrow\infty}{a_n} \ = \
x\left(\frac{R}{R-1}\right)^{R-1} \ \leq \ ex, \quad \mbox{and}
    \quad \limsup_{n\rightarrow\infty}{b_n} \ = \ 4R^Re^{-x},
    \label{eqn:anLimSup}
\end{equation}
by Corollary \ref{cor:patchLim}. Therefore by Lemma \ref{lem:KSV},
when $4R^Re^{-x}<1$, we have
\begin{eqnarray}
\nu^*(R) & \leq & \frac{ex}{1-4e^{-x}R^R}\nu^*(1).
    \nonumber 
\end{eqnarray}
Setting $f(x)=\frac{ex}{1-4e^{-x}R^R}$ and minimizing over $x>0$
such that $4e^{-x}R^R<1$, the derivative of $f$ is
$$ f'(x) \ = \ e\frac{\left(1-4(1+x)e^{-x}R^R\right)}
    {\left(1-4e^{-x}R^R\right)^2}.$$
The numerator $\left(1-4(1+x)e^{-x}R^R\right)$ has two roots, one
positive and one negative, and $f(x)$ reaches its minimum at the
positive root. Let this root be $x_0$, for which
$4e^{-x_0}R^R=\frac{1}{(1+x_0)}$, and so
\begin{eqnarray}
\nu^*(R) & \leq & e(x_0+1)\nu^*(1).  \nonumber
\end{eqnarray}
Since $\left(1-4(1+x)e^{-x}R^R\right)$ is negative on $[0,x_0)$
and increasing at $x_0$, we can bound $x_0$ slightly above by
choosing an approximation for $x_0$ which yields a positive value
in the numerator of $f'(x)$. Choosing
$x_0=(R\log{R}+\log{R}+\log{\log{R}}+3)$ ensures for $R\geq 2$
that $e^{x_0}>4(1+x_0)R^R$.  By Theorem 4.4.2 in \cite{CHLL97},
all optimal $(n,1)$-codes with length $n\geq 3$ are normal; and by
Theorem 12.4.11 in \cite{CHLL97}, $\mu^*(1)=1$; these results
allow the replacement of $\nu^*(1)$ with 1 to obtain the desired
result.
\end{proof}

\section{Asymptotic density of normal asymmetric codes}

We now present the asymmetric version of Theorem
\ref{thm:densitySym}, that is, a bound on the asymptotic
worst-case density of constant radius normal asymmetric codes. The
proof proceeds along the lines of that of the symmetric case, with
the most notable deviation occurring in the probabilistic
selection of the norm-patched asymmetric code $(S,T)$ due to a
more complicated definition of $T$. However, we obtain a
simplified asymptotic upper bound on $|T|$ which allows us to
employ the same analysis on the recursive ASDS construction as
before. The proof of Theorem \ref{thm:densityAsym} follows that of
Corollary \ref{cor:patchLimAsy}.

\begin{theorem}\label{thm:densityAsym}  Let $R\geq 2$.  Then
$$ \nu_+^*(R) \ \leq \ e(R\log{R}+\log{R}+\log{\log{R}}+4).
$$
\end{theorem}

Because of the asymmetry of the covering condition for
$(n,R)^+$-codes, we prefer to concentrate on the vast majority of
vertices of $Q_n$ which have weight close to $n/2$. Define a
vector $u\in Q_n$ to be {\em rare} if
$|w(u)-n/2|>\sqrt{2(R+1)n\ln{n}}$, and define
\begin{eqnarray}
hi(n,R) & := &
    \min{\left\{n,\left\lfloor\big(n+\sqrt{2(R+1)n\ln{n}}\big)/2
    \right\rfloor\right\}},\quad \mbox{and} \nonumber \\
lo(n,R) & := &
    \max{\left\{0,\left\lceil\big(n-\sqrt{2(R+1)n\ln{n}}\big)/2
\right\rceil\right\}}. \nonumber
\end{eqnarray}
Then the set of rare vectors of $Q_n$ (with respect to asymmetric
radius $R$) is
$$
Q^{rare}_n \ := \
    \left\{u\in Q_n: w(u)<lo(n,R) \mbox{ or } w(u)>hi(n,R) \right\}.
$$
The Chernoff bound states that the number of vertices $u\in Q_n$
with $w(u)>(n+\sqrt{j\cdot n\ln{n}})/2$ is at most $2^n
n^{-j^2/2}$ (cf.\ \cite[Theorem A.1.1]{AS00}). Thus
$|Q^{rare}_n|<2^{n+1}n^{-R-1}\in O(2^n n^{-R-1})$, which would
have density $O(1/n)$ as a $(n,R)^+$-code, except that for all but
finitely many $n$, $|Q^{rare}_n|$ doesn't downward $R$-cover
$Q_n$.

\begin{lemma}[Selection of a norm-patched asymmetric code]
\label{lem:delMethodAsy} For every positive constant $x$ and for
positive integers $N\leq n$, there exist (disjoint) sets
$S_0\subseteq \shc{(Q_n)}{0}{n}$ and $S_1\subseteq
\shc{(Q_n)}{1}{n}$ each of size at most
$$
    \frac{x2^{n-1}}
    {b_{n-1}^+\left(hi(n,R),\left\lfloor\frac{N-1}{2}\right\rfloor\right)
    +b_{n-1}^+\left(hi(n,R),\left\lceil\frac{N-1}{2}\right\rceil-1\right)},
$$
and a set $T\subseteq Q_n$ of size at most $\tau^+(n,N,x):=$
\begin{align}
O( & 2^n n^{-R-1}) +
    2^{n} \left[\sum_{i=0}^{N-1}
        \exp\left(-x\frac{b^+_{n-1}(hi(n,R),i-1)+b^+_{n-1}(hi(n,R),N-i-1)}
    {b_{n-1}^+\left(hi(n,R),\left\lfloor\frac{N-1}{2}\right\rfloor\right)
    +b_{n-1}^+\left(hi(n,R),\left\lceil\frac{N-1}{2}\right\rceil-1\right)}
    \right. \right. \nonumber\\
& \ + \ \left.\frac{b^+_{n-1}(hi(n,R),i-1)+b^+_{n-1}(hi(n,R),N-i-1)}
    {2^{n-1}}\right) \nonumber \\
& \ + \ \left. \exp\left(-x\frac{b^+_{n-1}(hi(n,R),N-1)}
    {b_{n-1}^+\left(hi(n,R),\left\lfloor\frac{N-1}{2}\right\rfloor\right)
    +b_{n-1}^+\left(hi(n,R),\left\lceil\frac{N-1}{2}\right\rceil-1\right)}
    \right.\right.\nonumber \\
& \ + \ \left.\left. \frac{b^+_{n-1}\left(hi(n,R)-1,\left\lfloor
    \frac{N-1}{2}\right\rfloor\right)}{2^{n-1}}
    \right)\right] \nonumber \\
& + \ 2^n
    \exp\left(-x\frac{b^+_{n-1}\left(hi(n,R)-1,\left\lfloor
    \frac{N-1}{2}\right\rfloor\right)}
    {b_{n-1}^+\left(hi(n,R),\left\lfloor\frac{N-1}{2}\right\rfloor\right)
    +b_{n-1}^+\left(hi(n,R),\left\lceil\frac{N-1}{2}\right\rceil-1\right)}
    \right.\nonumber\\
& \ + \ \left.\frac{b^+_{n-1}\left(hi(n,R)-1,\left\lfloor
    \frac{N-1}{2}\right\rfloor\right)}{2^{n-1}}\right)
\label{eqn:expPatchFullAsy}
\end{align}
such that $(S_0\cup S_1,T)$ is a balanced norm $N$-patched
asymmetric code.
\end{lemma}
\begin{proof}
Let
$$ k=\left\lfloor\frac{x2^{n-1}}
    {b_{n-1}^+\left(hi(n,R),\left\lfloor\frac{N-1}{2}\right\rfloor\right)
    +b_{n-1}^+\left(hi(n,R),\left\lceil\frac{N-1}{2}\right\rceil-1\right)}
    \right\rfloor,
$$
And choose uniformly at random subsets $S_0\subseteq
\shc{(Q_n)}{0}{n}$ and $S_1\subseteq \shc{(Q_n)}{1}{n}$ each of
size $k$. A vector $v\in \shc{(Q_n)}{0}{n}$ is missed by $S$ if
$d^+(v,S_0)+d^+(v,S_1)>N$; otherwise, there exists an
$i\in\{0,1,\ldots,N\}$ such that $d^+(v,S_0)=i$ and
$d^+(v,S_1)\leq N-i$.  A vector $v\in \shc{(Q_n)}{1}{n}$ is missed
by $S$ provided $2d^+(v,S_1)+1>N$. We classify the missed vertices
as follows:
\begin{eqnarray}
B^+_{0,1} & := & \{u\in\shc{(Q_n)}{0}{n}\setminus Q^{rare}_n
    :0\leq d^+(u,S_0)< N
    \mbox{ and } d^+(u,S_{1})>N-d^+(u,S_0)\} \nonumber \\
B^+_{0,2} & := & \{u\in\shc{(Q_n)}{0}{n}\setminus Q^{rare}_n
    :d^+(u,S_0)\geq N \} \nonumber \\
B^+_1 & := & \{u\in \shc{(Q_n)}{1}{n}\setminus Q^{rare}_n
    :2d^+(u,S_1)+1>N\}. \nonumber
\end{eqnarray}
Let the patch be the balanced set
$$ T \ = \ Q^{rare}_n\cup \left[\big(B^+_{0,1}\cup B^+_{0,2}\cup B^+_1\big)+
    (e_n\cup 0)\right];
$$
then $T$ contains all missed vectors, and $(S,T)$ is a balanced
norm $N$-patched asymmetric code.
By linearity of expectation,
\begin{eqnarray}
\E(|T|) & = & |Q^{rare}_n|+2\E(|B^+_{0,1}|)+2\E(|B^+_{0,2}|)+
    2\E(|B^+_{1}|)\nonumber \\
& = & O(2^n n^{-R-1}) +
    2\sum_{v\in\shc{(Q_n)}{0}{n}\setminus Q^{rare}_n}
    \sum_{i=0}^{N-1} \mbox{Pr}[d^+(v,S_0)=i]
    \cdot\mbox{Pr}[d^+(v,S_{1})>N-i]    \nonumber\\
& & + \ 2\sum_{v\in\shc{(Q_n)}{0}{n}\setminus Q^{rare}_n}
    \mbox{Pr}[d^+(v,S_0)\geq N]
    \ + \ 2\sum_{v\in\shc{(Q_n)}{1}{n}\setminus Q^{rare}_n}
    \mbox{Pr}[2d^+(v,S_1)+1 > N].
        \nonumber
\end{eqnarray}
Similar to the symmetric case, replacing $\mbox{Pr}[d^+(v,S_0)=i]$
above with $\mbox{Pr}[d^+(v,S_0)>i-1]$ yields a good upper bound
for $\E(|T|)$.
Using the definition of asymmetric distance, for any vector
$v\in\shc{(Q_n)}{0}{n}$ of weight $l$ and any $i$,
\begin{eqnarray}
\mbox{Pr}[d^+(v,S_0)>i] & = & \binom{2^{n-1}-b^+_{n-1}(l,i)}{k}
        \Big/ \binom{2^{n-1}}{k}, \ \mbox{and}\nonumber \\
\mbox{Pr}[d^+(v,S_1)>i] & = & \binom{2^{n-1}-b^+_{n-1}(l,i-1)}{k}
        \Big/ \binom{2^{n-1}}{k}. \nonumber
\end{eqnarray}
Similarly, for any $v\in\shc{(Q_n)}{1}{n}$ with weight $l$ and any $i$,
\begin{eqnarray}
\mbox{Pr}[d^+(v,S_1)>i] & = & \binom{2^{n-1}-b^+_{n-1}(l-1,i)}{k}
        \Big/ \binom{2^{n-1}}{k}. \nonumber 
\end{eqnarray}
Again using the estimate $\binom{m-d}{k}/\binom{m}{k}\leq
e^{-kd/m}$ in (\ref{eqn:binomEst}), this allows a regrouping of
the expression for $\E(|T|)$ by weight of $v$. We have
\begin{eqnarray}
\E(|T|) &  \leq & O(2^n n^{-R-1}) +
    2\sum_{l=lo(n,R)}^{hi(n,R)}
        \binom{n-1}{l} \left[\sum_{i=0}^{N-1}
        \exp\left(-k\frac{b^+_{n-1}(l,i-1)}{2^{n-1}}
        -k\frac{b^+_{n-1}(l,N-i-1)}{2^{n-1}}\right)
        \right. \nonumber\\
& & + \ \left. \exp\left(-k\frac{b^+_{n-1}(l,N-1)}{2^{n-1}}\right)\right]
\ + \ 2\sum_{l=lo(n,R)}^{hi(n,R)} \binom{n-1}{l-1}
    \exp\left(-k\frac{b^+_{n-1}\left(l-1,\left\lfloor
    \frac{N-1}{2}\right\rfloor\right)}{2^{n-1}}\right)  \nonumber\\
& \leq & O(2^n n^{-R-1}) +
    2^{n} \left[\sum_{i=0}^{N-1}
        \exp\left(-k\frac{b^+_{n-1}(hi(n,R),i-1)+b^+_{n-1}(hi(n,R),N-i-1)}
    {2^{n-1}}\right)
        \right. \nonumber\\
& & + \ \left. \exp\left(-k\frac{b^+_{n-1}(hi(n,R),N-1)}{2^{n-1}}\right)\right]
\ + \ 2^n
    \exp\left(-k\frac{b^+_{n-1}\left(hi(n,R)-1,\left\lfloor
    \frac{N-1}{2}\right\rfloor\right)}{2^{n-1}}\right)  \nonumber\\
& \leq & \tau^+(n,N,x). \nonumber
\end{eqnarray}
Since there exists a $T$ with size at most $\E(|T|)$, the result follows.
\end{proof}

Just as in the symmetric case, what is important about Lemma
\ref{lem:delMethodAsy} is the asymptotic behavior of $|T|$ as $n$
tends to infinity.  Accordingly, we have the following corollary.
\begin{cor}\label{cor:patchLimAsy}
Let $N\geq 2$ be fixed.
Then the asymptotic size of $S:=S_0\cup S_1$
in Lemma \ref{lem:delMethodAsy} is
$$
|S| \ \sim \
    \left\{\begin{array}{cl}
    x \frac{2^{n}}{(n/2)^R/R!}
    , & \mbox{if $N=2R+1$ is odd} \\
    x\frac{2^{n-1}}{(n/2)^{R-1}/(R-1)!}
    , & \mbox{if $N=2R$ is even,}
    \end{array}\right.
$$
and the size of the patch $T$ is bounded above asymptotically by
$$ \tau^+(n,N,x) \ \sim \
    \left\{\begin{array}{cl}
    3\cdot 2^{n}e^{-x}, & \mbox{if $N=2R+1$ is odd} \\
    2^{n}(e^{-x}+e^{-x/2}), & \mbox{if $N=2R$ is even.}
    \end{array}\right.
$$
\end{cor}
\begin{proof}
The asymptotic size of an upward asymmetric $R$-ball $B_n^+(v,R)$
for constant $R$ where $w(v)=hi(n,R)$ is
$$
    b_n^+(hi(n,R),R) \ = \ \binom{lo(n,R)}{\leq R} \ \sim \
    \frac{(n/2)^R}{R!}.
$$
The calculation for $|S|$ is now easily verified.
The proof of the bound for $T$  proceeds, similarly to the proof of
Cor.\ \ref{cor:patchLim}, by  identifying what exponential terms
$\exp(\cdot)$ in (\ref{eqn:expPatchFullAsy})
are not swallowed in the limit.
\end{proof}

\begin{proof}[Proof of Theorem \ref{thm:densityAsym}]
Let $n\geq R$, $n_1=\lfloor n/R \rfloor$, and $n_1'=n-n_1+1$ as in
the proof of Theorem \ref{thm:densitySym}.  We use Lemma
\ref{lem:delMethodAsy} to select a length $n_1'$ balanced norm
$(2R-1)$-patched asymmetric code $(S,T)$, where $|S|$ and $|T|$
are bounded above as given in the lemma.  Let $K_1$ be an optimal
normal $(n_1,1)^+$-code, and let $K_2$ be an optimal normal
$(n_1,R)^+$-code.  By Theorem \ref{thm:ASDS}, the ASDS of $(S,T)$
with $(K_1,K_2)$ has length $n$ and norm $2R+1$.  Therefore there
exists a normal $(n,R)^+$-code with size at most
$|(S,T)\dot\boxplus(K_1,K_2)|$, and so
\begin{align}
\nu_+( & n,R)  \ \leq \
    \left(\frac{|S||K_1|}{2}+\frac{|T||K_2|}{2}\right)\frac{\binom{\lfloor
    n/2\rfloor}{\leq R}}{2^n} \nonumber\\
& \leq \ \frac{1}{2}
    \frac{x2^{n_1'}}
    {b_{n_{1}'-1}^+\left(hi(n_1',R-1),R-1\right)
    +b_{n_1'-1}^+\left(hi(n_1',R-1),R-2\right)}
    \nu_+(n_1,1)
    \frac{2^{n_1}\binom{\lfloor
    n/2\rfloor}{\leq R}}{\binom{\lfloor n_1/2\rfloor}{\leq 1}2^n}
    \nonumber \\
& \ + \ \tau^+(n_1',2R-1,x)
        \frac{1}{2}\nu(n_1,R)
    \frac{2^{n_1}\binom{\lfloor
    n/2\rfloor}{\leq R}}{\binom{\lfloor n_1/2\rfloor }{\leq R}2^n}.
    \nonumber
\end{align}
Define $s_n:=\nu_+(n,R)$, $f_n:=\nu_+(n,1)$,
\begin{eqnarray} a_n & := &
\frac{1}{2}\frac{x2^{n_1'}}
    {b_{n_{1}'-1}^+\left(hi(n_1',R-1),R-1\right)
    +b_{n_1'-1}^+\left(hi(n_1',R-1),R-2\right)}
    \frac{2^{n_1}\binom{n/2}{\leq R}}{\binom{n_1/2}{\leq 1}2^n}, \quad
        \mbox{and}
        \nonumber\\
b_n & := & \frac{1}{2}\frac{2^{n_1}\binom{n/2}{\leq R}}
    {\binom{n_1/2}{\leq R}2^n}
        \tau^+(n_1',2R-1,x); \nonumber
\end{eqnarray}
note that
\begin{eqnarray}
\limsup_{n\rightarrow\infty}{a_n} & = & x
    \left(\frac{R}{R-1}\right)^{R-1} \ \leq \ ex, \quad \mbox{and}
    \quad \limsup_{n\rightarrow\infty}{b_n} \ = \ 3R^Re^{-x},
\nonumber
\end{eqnarray}
by Corollary \ref{cor:patchLimAsy}. Therefore by Lemma
\ref{lem:KSV}, when $3R^Re^{-x}<1$, we have
\begin{eqnarray}
\nu^*_+(R) & \leq & \frac{ex}{1-3e^{-x}R^R}\nu^*_+(1). \nonumber
\end{eqnarray}
Similar to the proof of Theorem \ref{thm:densitySym}, letting
$x_0=R\log{R}+\log{R}+\log\log{R}+3$ ensures for $R\geq 2$ that
the denominator of the right-hand side is positive, and gives the
desired result.
\end{proof}

\section{Open questions}
The primary open question, in the author's opinion, is the value
of $\mu^*_+(1)$, the asymptotic worst-case density of radius 1
asymmetric covering codes, for which we believe no respectable
upper bound has been published. This question is likely to be
quite hard, as it is related to the question of finding {\em
covering numbers}, specifically, the smallest number of
$l$-subsets of $[n]$ which contain all $(l-1)$-subsets of $[n]$
(cf.\ \cite{ARS03}). A more routine open question is to determine
for which values of $n$ and $R$ the ADS or ASDS constructions
yield best known upper bounds on $K^+(n,R)$. In general, the best
known lower and upper bounds on $K^+(n,R)$ (see \cite{CEK02,
ARS03, E03w}) are still open to significant improvement.

\section*{Acknowledgement}

Thanks are due to Iiro Honkala and Simon Litsyn for assistance in
identifying previous results.



\end{document}